\documentclass{amsart}
\usepackage{amssymb,graphics}
\usepackage{epstopdf}

\newcommand{\dontprint}[1]{\relax}
\newcommand{\bx}[1]
{\begin{picture}(8,8)\put(0.5,-1.0){\framebox(7,7){$\scriptstyle{#1}$}}\end{picture}}
\title[A few things I learnt from J\"urgen Moser]
{A few things I learnt from J\"urgen Moser}
\author{Alexander P. Veselov}
\address{Department of Mathematical Sciences,
Loughborough University, Loughborough LE11 3TU, UK  and Landau
Institute for Theoretical Physics, Moscow, Russia}
\email{A.P.Veselov@lboro.ac.uk}
\begin{document}
\maketitle
\vspace{0.5cm}
\rightline{\em To dear memory of J\"urgen Moser}
\bigskip

\begin{abstract}
A few remarks on integrable dynamical systems  inspired by discussions with J\"urgen Moser and by his work.
\end{abstract}

\section{Introduction}
 
Nine years have passed since the last time I spoke to J\"urgen Moser but it is still painful for me to think that this will never happen again. That time it was in his house in Schwerzenbach near Zurich, which I was passing by returning from a conference in July 1999. It was after the metastases had been discovered in his body and he already knew his fate. He started, seeing a question in my eyes, by saying that he most probably had only up to Christmas time to live. Then he looked at the 1958 German record of St Matthew Passion I had brought to him and changed the subject by telling me about his time at Princeton, of his friendship with John Nash and many other interesting details of his life around 1958. It was the usual Moser I had been lucky to know for about 10 years. I am telling this because it shows how strong this man was, which in combination with his dignity and wisdom makes him unique to everybody who knew him well enough.

I do not ever remember Moser arguing. It does not mean that he did not make clear his opinion, he actually always did this immediately in a very straightforward fashion but in such a way that would not offend even the most sensitive person. However he had never repeated his point afterwards in contrast, say, to me. On every occasion that I remember this happening, I was wrong and he was right. Once this was about the Euler rigid body system. I said that it seems to have been explicitly integrated in elliptic functions only in the second half of the XIX-th century (I think I mentioned Hermite in this relation). Moser did not think that it was that late but when I insisted he reacted in his usual way: ``Really ?" The next day I found on my desk a library volume of Jacobi's Gesammelte Werke opened at the first page of Jacobi's  letter (more than 50 pages long !) from March 1850 to the Acad\'emie des Sciences de Paris with full integration in elliptic functions not only of the Euler system but of the whole frame motion. There were other cases I had to understand my foolishness myself but Moser always left it with me to find out the truth and never returned to this point. I was not an exception of course: many of the visitors of FIM at ETH, Zurich during the time of his Directorate probably had similar experiences.

My happiest memories of Moser are, when after his busy days, he invites me to his office to talk maths. He is sitting in his chair with glasses on and from time to time making notes. Taking notes was one of his habits, which the participants of the Analysis seminar at ETH remember very well: he did this at every talk he attended. To talk to Moser was a real joy but sadly strict timekeeping was another habit of his, so to prolong these happy hours I sometimes accompanied him on his way to Stadelhofen railway station. 

This note, written for the special issue of  ``Regular and Chaotic Dynamics" dedicated to 80th anniversary of Moser, consists of a few remarks on integrable systems, which were inspired by discussions with J\"urgen Moser and by his work. I have presented them in several talks, but have never published them before. 

\section{Integrable modulations: could the adiabatic invariants be exact ?}

One of the problems I discussed with J\"urgen Moser in the early 1990s, when we  studied the integrable discretisations of the classical systems \cite{MV}, was the following.

Consider the family of integrable Hamiltonian systems with Hamiltonian $H(p,q|A)$ depending on the parameters $A=(a_1, \dots, a_N).$  One can think of the classical Neumann problem of the motion on the unit sphere with the Hamiltonian
$$H= \frac{1}{2}p^2 + \frac{1}{2}(a_1 q_1^2+\dots+ a_Nq_N^2)$$ as an example.
Let $F_1(p,q|A), \dots, F_N(p,q|A)$ be the corresponding $N$ independent integrals in involution.

Now suppose that we make the parameters time-dependent: $A=A(t)$ with $A(t)$ periodic with period $T.$ In general this of course destroys integrability of the system. I was looking for the special cases, which I called {\it integrable modulations}, when the values of the integrals $F_i(p(t),q(t)|A(t))$, where $(p(t),q(t))$ is an arbitrary solution of the modulated system
\begin{equation}
\label{mod}
\frac{dp}{dt}= - \frac{\partial H}{\partial q}(p,q|A(t)), \quad \frac{dq}{dt}= \frac{\partial H}{\partial p}(p,q|A(t)),\end{equation}
 are also $T$-periodic in time with the same period $T.$ This implies that the time $T$-shift along the trajectories of the modulated system gives an integrable symplectic map with the same integrals $F_1, \dots, F_N.$

Moser said that he would formulate the question as follows:
{\it when are the adiabatic invariants exact ?}
If we mean by the adiabatic invariants the action variables (see \cite{A}), then this implies that the modulated dynamics preserves the Liouville tori foliation of the original integrable system, so we have simply non-uniform and non-straight motion on each of these tori.

At that time I had an example of such integrable modulation for the rigid body system, which is the main subject of this section. But let me start with a simpler example of this phenomenon: the harmonic oscillator.\footnote{I am grateful to Anatoly Neishtadt for pointing out this example to me.} 

If we consider the Hamiltonian of the harmonic oscillator in the form
$$H=\frac{1}{2} (p^2+ \omega^2 q^2),$$ then it is easy to see that it has no non-trivial integrable modulations at all. However, if we consider the equivalent Hamiltonian system with
\begin{equation}
\label{1}
H=\frac{1}{2} \omega (p^2+q^2),
\end{equation}
 then {\it any} modulation $\omega=\omega(t)$ is integrable.

Indeed, the modulated system
$$\frac{dp}{dt}= - \omega q, \quad \frac{dq}{dt}= \omega p$$
clearly has integral $$I=\frac{1}{2}(p^2 + q^2)$$ (which is the action variable of the system), and by changing time $t \rightarrow \tau(t): \frac{d\tau}{dt}= \omega(t)$ we reduce it to the non-modulated form
$$\frac{dp}{d\tau}= - q, \quad \frac{dq}{d\tau}= p.$$
The modulated system has the solutions $$p=R \cos \tau(t), \, q(t) = R \sin \tau(t),$$
where the angle variable $\varphi=\tau(t)$ depends on $t$ in the way determined by $\omega(t).$
Note that the Hamiltonian (\ref{1})  (which is not an integral of motion anymore) changes periodically with time provided the modulation $\omega(t)$ is periodic.

This example shows that in general for a prescribed dependence on the parameters in the Hamiltonian we may not have integrable modulations at all. In particular, it seems that such modulations do not exist for the Neumann system or for the closely related Jacobi system of geodesics on ellipsoids \cite{M}.
That is why I found interesting that for the Euler rigid body system such modulation does exists even in the $N$-dimensional version, which we are going to see now.

After Arnold \cite{A} the motion of an $N$-dimensional rigid body is commonly identified with the geodesic flow for the  left-invariant metric on the Lie group $SO(N)$
$$\mathcal L = tr \frac{dX}{dt} J \frac{dX}{dt}^T.$$
Here $X \in SO(N),$ $J=J^T$ is some positive symmetric matrix and $Y^T$ denotes the transposition of matrix $Y.$ The corresponding angular momentum $M \in so(N)^*$ satisfies
the {\it Euler equation}
\begin{equation}
\label{eu}
\frac{dM}{dt} = [M,\Omega],
\end{equation}
where $\Omega= X^{-1}\frac{dX}{dt} \in so(N)$ is the angular velocity related to $M$ by
\begin{equation}
\label{mdef}
M=\Omega J + J \Omega.
\end{equation}
If we introduce the linear operator $A: so(N) \rightarrow so(N)^*$ by
$A(\Omega)=\Omega J + J \Omega,$ then the Hamiltonian of the Euler equation can be written as
$$H=\frac{1}{2}tr M^TA^{-1}(M).$$

The proof of integrability of the corresponding system was found by Manakov \cite{Man}.
It is based on the important observation that the Euler equation can be rewritten in the Lax form

\begin{equation}
\label{leu}
\frac{dL}{dt} = [L,P],
\end{equation}
where 
\begin{equation}
\label{lpeu}
L=M +  \lambda J^2, \,\, P=\Omega +  \lambda J
\end{equation}
and $\lambda$ is an arbitrary parameter (usually called spectral). Indeed, the relation (\ref{leu}) is equivalent to
$$\frac{dM}{dt} = [M,\Omega],\,\, [J^2,\Omega]+[M,J]=0, \,\, [J^2, J]=0.$$
The third relation is trivial, the second one is equivalent to (\ref{mdef}), so we end up with the original Euler equations.
The Lax equation (\ref{leu}) implies that the coefficients of the characteristic polynomial of $L$
$$P(\lambda,\mu)= \det (L-\mu I) = \det (M+\lambda J^2 -\mu I)$$
are integrals of the systems. One can check that this gives enough integrals in involution to claim the Liouville integrability of the system and that the dynamics is linearisable on the Prym variety of the
spectral curve $$P(\lambda,\mu)=0.$$

Let me show now that Manakov's trick works also when $J$ is not a constant but depends explicitly on $t$ in such a way that
$$J^2 = J_0^2 + f(t) I$$ with an arbitrary scalar function $f(t),$ constant symmetric matrix $J_0$ and the identity matrix $I.$
Indeed, the corresponding modulated Euler equation has the same form (\ref{eu}), but with 
$t$-dependent $$J = (J_0^2 + f(t) I)^{1/2}.$$ 
The same arguments as above show that it can be rewritten as
$$\frac{dM}{dt} = [M+\lambda J^2,\Omega+ \lambda J].$$
For our particular $J$ we have
$$[M+ \lambda J^2 ,\Omega+\lambda J]= [M+ \lambda J_0^2 + f(t) \lambda I, \Omega+\lambda J]=[M+ \lambda J_0^2,\Omega+\lambda J],$$
which implies that the modulated Euler system can be written in the Lax form (\ref{leu})
with 
$$L= M + \lambda J_0^2, \,\, P = \Omega + \lambda  J.$$ 

This means that the integrals of the non-perturbed system given by the coefficients of the characteristic polynomial $P_0(\lambda, \mu)=\det (M+\lambda J_0^2 - \mu I)$ are preserved by the modulated system, so
$$J = (J_0^2 + f(t)I)^{1/2}$$ is indeed an integrable modulation. 
The dynamics can be described explicitly as certain motion on the Prym variety of the corresponding spectral curve $P_0(\lambda, \mu)=0.$

Note that the Hamiltonian as well all other integrals of the initial system given by the coefficients of 
$P(\lambda, \mu)=\det (M+\lambda J^2 - \mu I)$ do depend on $t$:
$$P(\lambda, \mu)= \det (M+\lambda J_0^2+f(t) \lambda I -\mu I)= P_0(\lambda, \mu-f(t)\lambda)$$
and their $t$-dependence will be periodic if the modulation $f(t)$ was periodic.

The corresponding $T$-shift gives an integrable discretisation of the Euler rigid body system, which is worthy of further investigation. It depends on the choice of $f(t)$ and corresponds to certain shifts on the Prym varieties. In comparison with the discretisation from \cite{MV} its generating function (Lagrangian) is probably a very complicated function (cf. Bobenko-Lorbeer-Suris paper \cite{BS}, where a different discretisation of the Euler top was proposed).

\section{Integrability and analysis: is the cat-map integrable ?}

During one of my visits to Z\"urich, when Sergei Tabachnikov was also a guest of FIM, we thought that it would be a good idea to interview Moser for the Russian journal for high-school students ``Kvant" (American version ``Quantum"), of which Sergei was one of the Editors. Moser agreed without hesitation. We gave him around a dozen questions, which we thought he would simply write the answers to. Instead he invited us to come to talk about this. I still can not forgive myself that we were not ready for this and not able to prepare anything for publication afterwards. Because to say that this was interesting is to say nothing. 

One of the questions was about mathematics and where he sees himself. Moser was very reluctant to divide mathematics into separate disciplines and considered himself simply as a mathematician. However he added that now he is not as active as before and left for himself a small garden, where he tries to follow all the main events. This garden was called ``Analysis."

Integrable systems is a crossroad area, naturally connected to geometry and algebra. Analysis is usually connected to it only through the KAM-theory. J\"urgen Moser is amongst the very few real analysts, who made crucial contribution to the modern theory of integrability.
Many of us learnt this theory from his brilliant Bressanone's lectures \cite{M}.

I would like to talk now about one phenomenon in this theory which I have understood only recently.
This is about the difference between smooth integrability and analytic integrability.
One of the most famous Moser results is that KAM-theory works in smooth category as well as in analytic one, which was shown earlier by Kolmogorov and Arnold.
Here is a simple example, inspired by Bolsinov and Taimanov \cite{BT}, which demonstrates how big is actually the gap between these two categories.

Consider the following {\it cat-map} (or Anosov map)
$$A: T^2 \rightarrow T^2,$$ where
$T^2 = \mathbf R^2/ \mathbf Z^2$ is the two-dimensional torus and
$$A =\left(\begin{array}{cc} 2 & 1 \\ 1 &
1
\end{array}\right)$$ or any other hyperbolic 
$A=\left(\begin{array}{cc} a_{11} & a_{12} \\ a_{21} &
a_{22}
\end{array}\right) \in SL(2,{\mathbb Z}).$

The classical Anosov result says that any hyperbolic torus automorphism is structurally stable in $C^1$ topology. This is the canonical example of a chaotic system with positive topological entropy, which in this case can be shown to be equal to  $ \log \lambda,$ where $\lambda>1$ is the largest eigenvalue of $A.$ 

I argue that there are several reasons to consider this system as integrable.\footnote{To avoid complaints about abuse of terminology I would like to quote Birkhoff's classical book \cite{Birk}, whose 1966 edition was introduced by Moser: ``When, however, one attempts to formulate a precise definition of integrability, many possibilities appear, each with a certain intrinsic theoretic interest."}
First of all, this is a linear map and therefore can be integrated explicitly by means of linear algebra.
There is a problem with integrals though, because the dynamics is ergodic on the torus.
However I claim that the {\it natural extension of cat-map} to the cotangent bundle of the torus
$$\hat A: T^*T^2 \rightarrow T^*T^2$$ is integrable in the Liouville sense with {\it smooth} integrals.

Indeed, in the eigen-coordinates $u,v$ of $A$ the extended map has the form
$$\hat A: (u,v) \rightarrow (\lambda u, \lambda^{-1} v), \quad (p_u, p_v) \rightarrow (\lambda^{-1}p_u, \lambda  p_v),$$
so we have two integrals in involution:
$F_1 = p_u p_v$ and
$$F_2 = e^{-\frac{1}{p_u^2p_v^2}} \sin(2\pi \frac {\log p_u^2}{\log \lambda^2}).$$
Note that the second integral is smooth but not analytic.

The original cat-map is sitting at the {\it critical level} $p_u = p_v = 0.$
This is another lesson of this example (cf. Bolsinov-Taimanov \cite{BT}): Liouville-Arnold theorem
does not tell us anything about what happens at the degenerate level of the integrals. In smooth category one might find a fully developed chaos there !

There is of course another, {\it topological} aspect of this example. Topological obstructions to integrability have been discussed in the literature starting from the important work by Kozlov \cite{K}, but here I would like to make a different point. Let us look at the cat-map again. On the universal covering (which is the plane in the torus case) the dynamics is indeed integrable. Only after projection to the torus do we have ergodicity. A similar example is the geodesic flow on the surfaces of genus $g>1,$ which is chaotic on the surface but perfectly integrable on the hyperbolic plane covering it.
I would call such chaos {\it topological}. It has nothing to do with solvability.
There are many examples of exactly solvable systems, which are ergodic. 

Let me discuss one of them (probably the oldest one) due to Claude Gaspar Bachet (1581-1638), mathematician, poet, linguist, widely known by his translation of Diophantus' ``Arithmetica," where Fermat made his famous marginal notes. Looking for the rational solutions to the Diophantine equation
$$y^2 - x^3 = c,$$ Bachet observed that if $(x, y)$ is a solution, then so is
$$(x', y') = (\frac{x^4-8cx}{4y^2}, \frac{-x^6-20cx^3 + 8c^2}{8y^3}).$$
This allows to construct infinitely many solutions to this equation.
For example, for $c=-2$ we have
$$(3,5) \rightarrow (\frac{129}{100}, -\frac{383}{1000}) \rightarrow (\frac{2340922881}{7660^2}, \frac{113259286337292}{7660^3})\rightarrow\dots$$
The rational map $B: \mathbf CP^1 \rightarrow \mathbf CP^1$
$$x \rightarrow B(x) = \frac{x^4 -8cx}{4(x^3 + c)}$$
is called the {\it Bachet map.}
Geometrically this is the {\it duplication map} for the elliptic curve $E$ given by
$y^2 - x^3 =c.$
Indeed $$\wp(2z) = B(\wp(z)),$$
where $\wp(z)$ is the corresponding Weierstrass elliptic function,
so on $E$ we have the map $z \rightarrow 2z.$
The Bachet map is actually the most chaotic among the rational maps in the sense that
its Julia set coincides with the whole of $\mathbf CP^1.$
In complex dynamics this important example of chaotic map is usually ascribed to Latt\`es, who worked at the beginning of the last century (see e.g. Milnor \cite{Mil}). 

Again I would argue that this map has all the reasons to be considered as integrable.
One is exact solvability in terms of elliptic functions, which follows from the above arguments.
Another is the existence of infinitely many commuting rational maps (symmetries), which in the non-invertible maps replaces the integrals. Indeed, the rational maps $B_n$ defined by
$$\wp(n z) = B_n(\wp(z))$$ clearly commute with each other:
$$B_n B_m = B_m B_n, \quad B_2 = B.$$
Ritt showed that all such maps are related to some multiplication formula for elliptic functions. In many dimensions there are similar examples related to simple complex Lie algebras (see more details in \cite{V}).

\section{Geodesics on hyperboloids: Kn\"orrer map and geodesic equivalence as ``regularisations"}

Moser had a unique style as a lecturer. I attended one of his undergraduate courses on Dynamical Systems at ETH (initially to learn more German, but then just because I enjoyed it so much). Usually in the first hour he explained the result and the main ideas behind the proof using the simplest possible arguments and examples, then after the break full proofs with all the details were given. 
Not a single question was left unanswered and Moser, who usually kept time quite strictly, was always generous, making sure that that the matter was clarified completely. This combination of simplicity and clarity with mathematical honesty, not allowing him to hide the unwanted technicalities under the carpet, I found remarkable.

The same concerns his writing. Moser told me that the large number (333) of derivatives in his result on invariant curves in KAM theory was required just for the sake of simplicity of the arguments and he could easily reduce it (I think he said to 18) at the cost of extra technicalities. The paper \cite{M2} is a nice example of Moser's remarkable style, which is relevant to the problem I would like to discuss.

In the first part of this paper Moser considered Kepler's problem with Hamiltonian
$$H = \frac{1}{2} |p|^2 - \frac{1}{|q|}.$$ By changing the time variable $t$ to 
$$s = \int \frac{dt}{|q|}$$
and using stereographic projection, he showed that the energy surface $H=c$ with negative $c$ can be compactified  and that the regularised flow is equivalent simply to the geodesic flow on a unit sphere (three-dimensional in the usual case, when $q \in \mathbf R^3$). We note that because of the collision orbits (corresponding to $|q| \rightarrow 0$) the regularisation is really needed.  The change of time allows to continue the motion after the collision in a smooth way. Moser mentioned that such regularisation was essentially known to Levi-Civita and Sundman, but it looks like that this topological picture was clearly explained first in his paper. \footnote{On the quantum level this corresponds to the $SO(4)$ symmetry of the Kepler problem with negative energy discovered by V. Fock in 1935.}

A simple remark I would like to make is related to the geodesic problem on hyperboloids.
All the algebraic manipulations are the same as in the well-studied Jacobi geodesic problem on ellipsoids. However the behaviour of geodesics in these two problems are quite different: in the Jacobi case the dynamics is quasi-periodic while on hyperboloids the geodesics in general go to infinity. 
I would like to explain here that this difference disappears if one makes a ``compactification-regularisation" of the geodesic problem on hyperboloid using Kn\"orrer's transformation \cite{Kn}. 

Let me recall first Kn\"orrer's result (which I also learnt from Moser \cite{M3}). Let 
\begin{equation}
\label{Q}
(A^{-1}x,x)=1
\end{equation}
be the equation of a quadric $Q$ in the $(n+1)$-dimensional Euclidean vector space $\mathbb R^{n+1}$, which for the beginning will be assumed to be an ellipsoid. Geodesics on $Q$ in the natural (length) parameter $s$ are described by the equation
$$x'' = \lambda Bx,$$ where $'$ denotes the derivative with respect to $s$, $B= A^{-1}$, $\lambda$ is the Lagrange multiplier determined by the constraint (\ref{Q}):
$$(Bx, x')=0, \quad (Bx, x'')+(Bx',x')=0, \quad \lambda (Bx, Bx)+(Bx',x')=0,$$ 
\begin{equation}
\label{lambda}
\lambda=-\frac{(Bx',x')}{|Bx|^2}.
\end{equation}
Following Kn\"orrer \cite{Kn}, consider the change of parametrisation
\begin{equation}
\label{tau}
s \rightarrow \tau = \int \alpha(s)ds,
\end{equation}
where 
$\alpha^2 = -\lambda.$ 
Kn\"orrer observed the remarkable fact that the normal vector to the geodesic 
\begin{equation}
\label{G}
q = \frac{Bx}{|Bx|}
\end{equation}
 in the new time $\tau$ 
satisfies the classical Neumann system:
\begin{equation}
\label{N}
\frac{d^{2}q}{d\tau^2}= -Bq +\mu q,\,\, q \in\mathbb R^{n+1},\,\, |q|=1.
\end{equation}
In other words, the Gauss map (\ref{G}) in combination with re-parametrisation (\ref{tau}) transforms the geodesics into the trajectories of the Neumann system on the unit sphere in the harmonic field with the potential
$$U(q)=\frac{1}{2}(Bq,q).$$
The corresponding orbits have one of the integrals fixed: in Moser's notations \cite{M3}
\begin{equation}
\label{rel}
\Psi_0(\frac{dq}{dt},q)=0,
\end{equation}
where $$\Psi_0(u,v)=(1-(Au,u)) A(v,v) + (Au,v)^2, A=B^{-1}.$$
Thus, Kn\"orrer transformation gives a bijection between the unit (co)tangent bundle $\mathcal M$ to ellipsoid $Q$
and the subvariety $\mathcal N$ of the (co)tangent bundle to the unit sphere given by relation (\ref{rel}),
mapping the geodesics to Neumann orbits.

Let us look now at the case when the quadric $Q$ is a hyperboloid. In that case generically the geodesics will go to infinity and the corresponding $\lambda$ tends to 0 as $s \rightarrow \infty.$
To see the rate of decay we note that the geodesic problem has the following Joachimsthal integral
$$F=|Bx|^2(Bx',x').$$ Indeed,
$$F'=2(Bx',Bx)(Bx',x')+2|Bx|^2(Bx',x'')$$
$$=2(Bx',Bx)(Bx',x')+2\lambda |Bx|^2(Bx',Bx)=2(Bx',Bx)((Bx',x')+\lambda |Bx|^2) \equiv 0.$$
Using this integral we can write 
$$\lambda = -\frac{(Bx',x')}{|Bx|^2}=-\frac{F}{|Bx|^4}.$$
Because at infinity $|Bx| \approx cs$ we see that $\lambda \approx C s^{-4}$ and thus
$\alpha = \sqrt -\lambda$ decays as $s^{-2}$ (one should be careful with the sign of $\lambda$ in the hyperboloid case, but this can be done in all the cases). This means that the integral $\tau = \int \alpha(s)ds$
is convergent at infinity, so it takes {\it finite time} in the new time variable $\tau$ to reach infinity.
The Neumann system tells us what happens afterwards: the geodesic will return from the "same infinity" (in general, along a different path), then go to a different one and so on. Note that the Neumann system does not change when $B \rightarrow B-zI$ for any $z$, so it does not feel the change from ellipsoid to hyperboloid, which affects only the value of integral $\Psi_0$ (cf. Moser's comments in section 3.5 in \cite{M3}).

The Kn\"orrer transformation makes the necessary compactification and defines a smooth dynamics there. To see what happens geometrically, it is convenient to use the dual quadric $Q^*$ given by
$$(Ay,y)=(B^{-1}y,y)=1,\,\, y = Bx \in\mathbb R^{n+1}.$$ Then the Gauss map becomes simply the {\it central projection} of the hyperboloid $Q^*$ to a unit sphere:
$$p: y \rightarrow q=\frac{y}{|y|}.$$ The image is an open set with the boundary $\mathcal B$ given by the intersection of the sphere with the asymptotic cone 
$(Aq,q)=0.$ The points at the boundary correspond to the "infinities" of the hyperboloid. Note that when $(Aq,q)=0$  the relation (\ref{rel}) reduces to $$(Aq, \frac{dq}{dt})=0,$$ which means that the velocity vector must be tangent to the boundary.

We can see this already in the simplest example of hyperbola 
$$b_0x_0^2 + b_1x_1^2=1$$
with $b_0>0,\,\, b_1<0.$ In this case the regularised flow will be simply periodic back and forth motion on one branch of the hyperbola. The unit tangent bundle $\mathcal M$ in this case is topologically equivalent to the disjoint union of two open intervals, while on the Neumann side we have a circle, which is a compactification of $\mathcal M$ with two extra points, corresponding to two "infinities" of the hyperbola.

The two-dimensional case is more interesting. In particular, in the case of a one-sheeted hyperboloid we have
two families of straight lines, which are the simplest geodesics. For them (and only for them) we have $(Bx',x')=0$ and correspondingly $\lambda \equiv 0,$ which means that the Kn\"orrer transformation formally can not be applied. The central projection maps them into large semicircles, corresponding to the large energy limit of the Neumann system. This shows that we may have something extra on the geodesic side as well. 

There is another, more geometrical way to compactify the geodesic system on hyperboloid, which I have realised after discussion of this problem with Alexey Bolsinov. It is based on the notion of geodesic equivalence of the metrics. Two metrics $ds_1^2$ and $ds_2^2$ (Riemannian or pseudo-Riemannian ) on a manifold $M^n$ are called {\it geodesically equivalent} if they have the same geodesics considered as {\it unparametrised} curves. A trivial example is given by the proportional metrics $ds_2^2 = c ds_1^2,$ but there are other examples as well. The problem of description of projectively equivalent metrics has a long history, going back to the work of Beltrami and Dini in XIX-th century (see the recent review \cite{TM}, where the role of integrability in this problem is emphasized).

We need the following example of two geodesically equivalent metrics, which, in spite of its classical look, seems to be discovered only recently independently by Tabachnikov \cite{T} and Matveev-Topalov \cite{TM}. Let $ds_1^2$ be the restriction of the Euclidean metric $$ds^2=\sum_{i=1}^{n+1} dx_i^2 = (dx, dx)$$ in $\mathbb R^{n+1}$ to the quadric $Q$ given by (\ref{Q}) (which is a standard metric on $Q$) and 
$ds_2^2$ be the restriction of the metric
\begin{equation}
\label{dr}
dr^2= \frac{b_1 dx_1^2 + b_2 dx_2^2 + \dots + b_{n+1} dx_{n+1}^2}{b_1^2x_1^2 + b_2^2 x_2^2 \dots + b_{n+1}^2 x_{n+1}^2}=\frac{1}{|Bx|^2}(Bdx,dx)
\end{equation}
to $Q.$ Then the claim is that these two metrics are geodesically equivalent (see \cite{T, TM}). This result was originally formulated for the ellipsoids, but it obviously holds for the hyperboloids as well although the second metric may become non-positive (cf. Theorem 4.13 in Khesin-Tabachnikov \cite{KT}).

The observation, which may possibly be new, is that in the hyperboloid case the second metric $ds_2^2$ can be naturally extended to a {\it regular metric on the projective closure} $\bar Q \subset \mathbb RP^{n+1}.$ The easiest way to see this is by direct calculation. In the projective coordinates $z_0:z_1:\dots :z_{n+1}$ the quadric $\bar Q$ is given by the equation
$$b_1z_1^2 + \dots + b_{n+1}z_{n+1}^2 = z_0^2,$$ where the original quadric $Q$ is the affine part in the chart, where $z_0 \neq 0$ and $x_k=z_k/z_0.$ In a different chart, where say $z_1 \neq 0,$
we can use the affine coordinates $y_k = z_k/z_1, \, k =0,\, 2,\,\dots,\, n+1$ such that
$$x_1=\frac{1}{y_0}, \,x_2=\frac{y_2}{y_0}, \dots, \, x_{n+1}=\frac{y_{n+1}}{y_0},$$
in which the equation of $\bar Q$ becomes 
$$b_1 + b_2 y_2^2 \dots + b_{n+1}y_{n+1}^2 = y_0^2.$$ 
One can easily check that the second metric $ds_2^2$ in these coordinates coincides with the restriction  to $\bar Q$ of the metric
\begin{equation}
\label{drtil}
d\tilde r^2 = \frac{-dy_0^2 + b_2 dy_2^2 + \dots + b_{n+1} dy_{n+1}^2}{b_1^2 + b_2^2 y_2^2 \dots + b_{n+1}^2 y_{n+1}^2}
\end{equation}
(note that outside $\bar Q$ the metrics $d\tilde r^2$ and $dr^2$ are different). A remarkable fact is that this metric is regular when $y_0=0,$ which means that $ds_2^2$ indeed can be regularly extended to the whole $\bar Q.$ 

When $n=2$ we have one-sheeted and  two-sheeted hyperboloids when $b_1 < 0 < b_2 < b_3$ and $b_1 < b_2 < 0 < b_3$ respectively. In the first case $\bar Q = T^2$ topologically is two-dimensional torus and the metric $ds_2^2$ has the signature $(1,1)$, in the second case $\bar Q = S^2$ is a topological sphere with the Riemannian metric $(-ds_2^2)$.

The relationship of the geodesic equivalence with Kn\"orrer transformation and the general geometrical picture is worthy of further study. We only mention here that the comparison of the formulas (\ref{lambda}), (\ref{tau}) with (\ref{dr}) shows that the corresponding reparametrisations are the same. 
It may give also an alternative explanation of the special role of metric (\ref{dr}) on the quadrics.

Another interesting question is to look at the properties of the spectral problem 
$$-\psi'' + E \lambda(s) \psi =0,$$ where $\lambda$ is the corresponding Lagrange multiplier (\ref{lambda})
(see \cite{Cao, V2}, where the ellipsoid case was considered).

\section*{Acknowledgements} 
I am very grateful to my Loughborough colleagues E. Ferapontov, I. Marshall,  A. Neishtadt and especially to A. Bolsinov for very fruitful discussions. I am also grateful to P. Grinevich and S. Abenda, who attracted my attention to the paper \cite{Cao}, and to M. Levi and S. Tabachnikov for supportive comments.

 This work has been partially supported by the
European Union through the FP6 Marie Curie RTN ENIGMA (Contract
number MRTN-CT-2004-5652) and ESF programme MISGAM and by the EPSRC (grant EP/E004008/1).

\end{document}